\documentclass{IEEEtran}
\usepackage{cite}
\usepackage{amsmath,amssymb,amsfonts}
\usepackage{algorithmic}
\usepackage{graphicx}
\usepackage{textcomp}
\usepackage{amsmath, bm}
\usepackage{cases}
\newtheorem{rem}{Remark}
\newtheorem{assumption}{Assumption}
\newtheorem{problem}{Problem}

\newtheorem{theorem}{Theorem}
\newtheorem{lemma}{Lemma}
\def\BibTeX{{\rm B\kern-.05em{\sc i\kern-.025em b}\kern-.08em
    T\kern-.1667em\lower.7ex\hbox{E}\kern-.125emX}}
\begin{document}
\title{Mismatched Disturbance Rejection Control for Second-Order Discrete-Time Systems}
\author{Shichao Lv, Kai Peng, Hongxia Wang, and Huanshui Zhang, \IEEEmembership{Senior Member, IEEE}
\thanks{Manuscript received Month xx, 2xxx; revised Month xx, xxxx; accepted Month x, xxxx. This work was supported by the Foundation for Innovative Research Groups of the National Natural Science Foundation of China (61821004), Major Basic Research of Natural Science Foundation of Shandong Province (ZR2021ZD14), High-level Talent Team Project of Qingdao West Coast New Area (RCTD-JC-2019-05), Key Research and Development Program of Shandong Province (2020CXGC01208), and Science and Technology Project of Qingdao West Coast New Area (2019-32,2020-20, 2020-1-4).}
\thanks{Shichao Lv, Hongxia Wang, and Huanshui Zhang are with the College of Electrical Engineering and Automation, Shandong University of Science and Technology, Qingdao, 266590, China (e-mail: lv\_sc2020@sdust.edu.cn; whx1123@163.com; hszhang@sdu.edu.cn). }

\thanks{Kai Peng is with the School of Power and Energy, Northwestern Polytechnical University, Xi{'}an, 710072, China (e-mail: pengkai@nwpu.edu.cn).}}

\maketitle

\begin{abstract}
This paper is concerned with mismatched disturbance rejection control for the second-order discrete-time systems. Different from previous work, the controllability of the system is applied to design the disturbance compensation gain, which does not require any coordinate transformations. Via this new idea, it is shown that disturbance in the regulated output is immediately and directly compensated in the case that the disturbance is known. When the disturbance is unknown, an extra generalized extended state observer is applied to design the controller. Two examples are given to show the effectiveness of the proposed methods. Numerical simulation shows that the designed controller has excellent disturbance rejection effect when the disturbance is known. The example with respect to the permanent-magnet direct current motor illustrates that the proposed control method for unknown disturbance rejection is effective.
\end{abstract}


\section{Introduction}
\label{sec:introduction}
Disturbances always bring negative influence on systems performance in practice. Disturbance is one of the most important factors influencing system performance. It thus is necessary for high-precision control to investigate the disturbance rejection technique. During the past four decades, some elegant disturbance attenuation techniques such as disturbance-observer-based control (DOBC) \cite{1372532}, extended state observer-based control (ESOBC) \cite{han1995extended}, disturbance accommodation control (DAC) \cite{johnson1971accomodation}, and composite hierarchical anti-disturbance control disturbance regulation control (CHADC) \cite{guo2014anti}, and so on, have been proposed and applied in various fields, see \cite{5572931,MADONSKI201518,9105115,roman2021hybrid,9345466,ding2021extended}.

The existing results mainly focus on treating matching disturbances \cite{johnson1970further,4391082,han2009pid,li2014disturbance}. Such disturbances always get involved the system via the same channel with the control input, or the influence of the disturbances can be equivalently transformed into the input channels in some way. Due to the aforementioned feature, the influence of the matching disturbances on the system can be counteracted directly or indirectly by appropriately designing the controller in the same channel (before or after transformed) as the disturbances are known. In the case the disturbances are unknown, disturbance observer (DOB)/extended state observer (ESO) was then applied to disturbance rejection \cite{5572931,7265050,2012Generalized}. Hence, the matching disturbance rejection technique has become relatively mature.

In contrast, mismatched disturbance rejection is more challenging. Mismatched disturbances extensively exist in the real world. Many practical systems like the permanent magnet synchronous motors, roll autopilots for missiles and flight control systems, are affected by the mismatched disturbances \cite{mohamed2007design,chwa2004compensation,chen2003nonlinear}. Essentially different from the matching disturbances, these disturbances act on the system via the different channel from the control input. Consequently, no matter what kind of control scheme is employed, it may be impossible to completely eliminate the influence of the mismatched disturbances on the system state \cite{isidori1985nonlinear}. Therefore, a practical consideration goals at eliminating the effects of mismatched disturbances from some variables of interests, i.e., from the output channel.

There are a number of methods in dealing with the mismatched disturbances \cite{2012Generalized,castillo2018enhanced,ginoya2013sliding,chen2016adrc,yang2012nonlinear,6129407,9339876}. In \cite{2012Generalized}, a method based on a generalized ESO is proposed for linear systems and it can eliminate the mismatched disturbances from the controlled output in steady state. Similar to \cite{2012Generalized}, \cite{castillo2018enhanced} weakens the restriction of disturbance and improves the disturbance reject effect by introducing the high-order derivative of the disturbance. For the nonlinear system with mismatched disturbance, a novel slide-mode control based on a generalized DOB is provided in \cite{ginoya2013sliding,6129407} and it can reject the mismatched disturbance in the output in steady state. Different from \cite{ginoya2013sliding,6129407}, \cite{yang2012nonlinear} treats the mismatched disturbance in the multi-input multi-output system with arbitrary disturbance relative degree. However, it should be noted that the mismatched disturbance could not be eliminated immediately and completely in the literatures like the case of matched disturbance. How to compensate the disturbance immediately and directly in the case of mismatching is remained to be solved.

This paper aims to proposing a controller design method to eliminate the mismatched disturbances in the output directly and completely. Different from those existing ideas, the controllability of the system is utilized to establish a connection between the control and the disturbance. Based on this connection, the disturbance rejection controller is designed. As the disturbances are known, the proposed controller can eliminate the mismatched disturbance in the output completely. Given that the disturbance is usually unknown, a modified controller is also presented based on our previous controller and a generalized ESO. The modified controller can remove the disturbance in the steady-state regulated output. It should be highlighted that the residual disturbance in the regulated output under the modified controller is induced by the error of the observer. In other words, the modified controller can eliminate the disturbance in the regulated output completely if the observer can provide the exact disturbance value at every instant.

The remainder of the paper is organized as follows.  Section \ref{sec2} introduces the problem of mismatched disturbance rejection control for second-order discrete-time systems. Method of mismatched known disturbance rejection control is proposed in Section \ref{sec3}. In Section \ref{sec4}, we modify the controller proposed in Section \ref{sec3} by introducing a generalized ESO to propose a controller for unknown disturbances. In Section \ref{sec5}, two examples are given to show the effectiveness of the proposed methods. The conclusions are drafted in Section \ref{sec6}.

\section{Problem Statement}\label{sec2}
\indent We consider the second-order single-input single-output (SISO) system with mismatched disturbances
\begin{numcases}{}\label{f2.1}
{\bm{{x}}(k+1)=\bm{Ax}(k)+\bm{b}_u u(k) +\bm{b}_d d(k)},\notag \\
{\bm{y}(k)=\bm{C}_m \bm{x}(k)},\\
{{y}_{o}(k)=\bm{c}_o \bm{x}(k)}\notag
\end{numcases}
where $\bm{x}(k)\in\bm{R}^2$, $u(k)\in \bm{R}$, $d(k)\in \bm{R}$, $\bm{y}(k)\in \bm{R}^r$, and ${y}_{o}(k)\in \bm{R}$ are the state, control input, disturbance, measurable output, and regulated output, respectively. $ \bm{A}$, $\bm{b}_u$, $\bm{b}_d$, $\bm{C}_m$, and $\bm{c}_o$ are given matrices with dimension $2\times 2$, $2\times 1$, $2\times1$, $r\times 2$, and $1\times 2$, respectively.
\begin{rem}
In (\ref{f2.1}), $d(k)$ represents known or unknown disturbance. The known disturbance includes the measurable disturbance and the disturbance models are available\cite{yang1994disturbance,10.1115/1.2896180}, etc. The unknown disturbance possibly includes external disturbances, un-modeled dynamics, parameter uncertainty, complex nonlinear dynamics, discretization error, and so on.
\end{rem}

\begin{rem}\label{rem:2}
The mismatched disturbance which implies that the disturbance $d(k)$ enters the system in a different channel from the control input $u(k) $\cite{6129407}. For the second-order system (\ref{f2.1}), that $d(k)$ is the mismatched disturbance means ${\bm{b}_u}'{\bm{b}_d}=0$.
\begin{assumption}\label{assumption.1}
${\bm{c}_o}{\bm{b}_u}=0$, which means we focus on removing the disturbance in the channel with the disturbance rather than with the control input.
\end{assumption}

\end{rem}

We will deal with the following two problems in this paper.
\begin{problem}\label{problem1}
For the system (\ref{f2.1}), assume the disturbance $d(k)$ is bounded and known. Seek $u(k)$ to stabilize the system and eliminate the disturbance of the regulated output immediately and completely.
\end{problem}

\begin{problem}\label{problem2}
For the system  (\ref{f2.1}), assume the disturbance $d(k)$ is unknown. Seek $u(k)$ to stabilize the system and remove the disturbance of the regulated output in steady state.

\end{problem}

\section{Known Disturbances Rejection Controller}\label{sec3}
\subsection{Control Law Design}

In order to stabilize the system and eliminate the disturbance of the regulated output simultaneously, the design of the controller can be divided as
\begin{align}\label{max.1}
u(k)={u}_{x}(k)+{u_{d}(k)}
\end{align}
where ${u}_{x}(k)$ and ${u_{d}(k)}$ are responsible for stabilizing the system (\ref{f2.1}) and compensating the disturbance of the regulated output, respectively.
${u}_{x}(k)$ and ${u_{d}(k)}$ will be given in the following theorem.

\begin{theorem}\label{the:1}
Suppose Assumption \ref{assumption.1} holds, the state is available, $(\bm{A}, \bm{b}_u)$ is controllable, and the disturbance is bounded and known.\\
1) Let $u_x(k)=\bm{Kx}(k)$, where the state feedback gain $\bm{K}$ is selected to guarantee that matrix $\bm{A}+\bm{b}_u\bm{K}$ is Schur;\\
2) Let
\begin{align}\label{3.4}
{u_{d}(k)}=-\bm{I}_1\bm{K}_p\bm{b}_d d(k)-\bm{I}_2\bm{K}_p\bm{b}_d d(k+1)
\end{align}
where $\bm{K}_p=\begin{bmatrix}{\bm{b}_u},(\bm{A}+\bm{b}_u\bm{K}){\bm{b}_u}\end{bmatrix}^{-1}$, $\bm{I}_1=\begin{bmatrix}1&0\end{bmatrix}$, $\bm{I}_2=\begin{bmatrix}0&1\end{bmatrix}$.\\
Under the composite control law (\ref{max.1}) with ${u_{x}(k)}$ in 1) and ${u_{d}(k)}$ in 2), there hold
\begin{itemize}
\item The closed-loop system (\ref{f2.1}) is bounded-input-bounded-state (BIBS) stable;
\item The known disturbance can be immediately and completely eliminated from the regulated output.
\end{itemize}

\end{theorem}

\begin{rem}
In the case  the disturbance is time invariant or slow time-varying, we may let $d(k+1)=d(k)$, then (\ref{3.4}) becomes:
\begin{equation}\label{k_d000}
\begin{aligned}
{u_{d}(k)}&=-\bm{I}_1\bm{K}_p\bm{b}_d d(k)-\bm{I}_2\bm{K}_p\bm{b}_d d(k)\notag\\
&=-(\bm{I}_1+\bm{I}_2)\bm{K}_p\bm{b}_d  d(k).
\end{aligned}
 \end{equation}

\end{rem}

\begin{IEEEproof}\label{theorem}
Firstly, we will show $\begin{bmatrix}{\bm{b}_u},(\bm{A}+\bm{b}_u\bm{K}){\bm{b}_u}\end{bmatrix}$ is nonsingular.
Since $(\bm{A},\bm{b}_u)$ is controllable, $(\bm{A}+\bm{b}_u\bm{K},\bm{b}_u)$ is controllable for any compatible $\bm{K}$. Hence, $\begin{bmatrix}{\bm{b}_u},(\bm{A}+\bm{b}_u\bm{K}){\bm{b}_u}\end{bmatrix}$ is full row rank. For the second-order SISO system (\ref{f2.1}), the full row rank of $\begin{bmatrix}{\bm{b}_u},(\bm{A}+\bm{b}_u\bm{K}){\bm{b}_u}\end{bmatrix}$ is equivalent its nonsingularity.

Secondly, we will show that the system (\ref{f2.1}) is BIBS stable under the proposed control law in Theorem \ref{the:1}.
Inserting this control law into the state equation in the system (\ref{f2.1}), leads to
\begin{align}\label{f2.01}
{{\bm{x}}}(k+1)=&\bm{A}\bm{x}(k)+\bm{b}_u({u}_{x}(k)+{u_{d}(k)})+\bm{b}_dd(k)\notag\\
=&(\bm{A}+\bm{b}_u\bm{K}){\bm{x}}(k)+\bm{b}_u(-\bm{I}_1\bm{K}_p\bm{b}_d d(k)\\
&-\bm{I}_2\bm{K}_p\bm{b}_d d(k+1))+\bm{b}_dd(k)\notag.
\end{align}
Since $d(k)$ is bounded. $\bm{b}_u(-\bm{I}_1\bm{K}_p\bm{b}_d d(k)-\bm{I}_2\bm{K}_p\bm{b}_d d(k+1))+\bm{b}_dd(k)$ is bounded. Note that $\bm{A}+\bm{b}_u\bm{K}$ is Schur matrix. According to Proposition 5.8 in \cite{d2013control}, the closed-loop system (\ref{f2.01}) is BIBS stable.

Thirdly, we will prove that $y_o(j)$ for $j\geq h$ in (\ref{f2.1}) is unaffected by disturbance $d(l)$ ($l\geq h$) under the disturbance rejection law proposed in Theorem \ref{the:1} if the disturbance is present from $h$ ($h>0$) instant.

In the case of $j=h$, there is no disturbance $d(h)$ in $y_o(h+1)$, which stems from the derivation below.
Plugging the control law given in Theorem \ref{the:1} into the state equation in the system (\ref{f2.1}), the state $\bm{x}(h)$ is expressed as
\begin{equation}\label{f2.09}
\begin{aligned}
{{\bm{x}}}(h)=&\bm{A}\bm{x}(h-1)+\bm{b}_u({u}_{x}(h-1)+{u_{d}(h-1)})\\
=&(\bm{A}+\bm{b}_u\bm{K}){\bm{x}}(h-1)-\bm{b}_u(\bm{I}_1\bm{K}_p\bm{b}_d d(k)\notag\\
&+\bm{I}_2\bm{K}_p\bm{b}_d d(k+1)),
\end{aligned}
 \end{equation}

where we have used the fact that the disturbance does not appear before $h$ instant.
Inserting the above equation into the regulated output equation in the system (\ref{f2.1}) and using ${\bm{c}_o}{\bm{b}_u}=0$ yield
\begin{equation}\label{f4.9}
\begin{aligned}
y_o(h)=&\bm{c}_o{{\bm{x}}}(h)=\bm{c}_o((\bm{A}+\bm{b}_u\bm{K}){\bm{x}}(h-1)\\
&-\bm{b}_u(\bm{I}_1\bm{K}_p\bm{b}_d d(k)+\bm{I}_2\bm{K}_p\bm{b}_d d(k+1)))\\
=&\bm{c}_o(\bm{A}+\bm{b}_u\bm{K}){\bm{x}}(h-1)-\bm{c}_o\bm{b}_u(\bm{I}_1\bm{K}_p\bm{b}_d d(k)\notag\\
&+\bm{I}_2\bm{K}_p\bm{b}_d d(k+1))\\
=&\bm{c}_o(\bm{A}+\bm{b}_u\bm{K}){\bm{x}}(h-1).
\end{aligned}
 \end{equation}
Obviously, there is no disturbance in $y_o(h)$.

In the case of $j = h+1$, it is easy to know
\begin{equation}\label{f2.8}
\begin{aligned}
{{\bm{x}}}(h+1)=&(\bm{A}+\bm{b}_u\bm{K}){\bm{x}}(h)+\bm{b}_u{u}_{d}(h) +\bm{b}_dd(h)\\
=&(\bm{A}+\bm{b}_u\bm{K})((\bm{A}+\bm{b}_u\bm{K}){\bm{x}}(h-1)\\
&+\bm{b}_u{u}_{d}(h-1))+\bm{b}_u{u}_{d}(h) +\bm{b}_dd(h)\\
=&(\bm{A}+\bm{b}_u\bm{K})^{2}{\bm{x}}(h-1)+\bm{b}_u{u}_{d}(h) \\
&+(\bm{A}+\bm{b}_u\bm{K})\bm{b}_u{u}_{d}(h-1)+\bm{b}_dd(h)\\\
=&(\bm{A}+\bm{b}_u\bm{K})^{2}{\bm{x}}(h-1)+\bm{b}_dd(h)\\
&+\begin{bmatrix}{\bm{b}_u},(\bm{A}+\bm{b}_u\bm{K}){\bm{b}_u}\cr\end{bmatrix}\begin{bmatrix}{u}_d(h)\cr{u}_{d}(h-1)\end{bmatrix}.
\end{aligned}
\end{equation}
Let $\bm{P}=\begin{bmatrix}{\bm{b}_u},(\bm{A}+\bm{b}_u\bm{K}){\bm{b}_u}\cr\end{bmatrix}$, substituting (\ref{3.4}) with $k=h-1$ and $k=h$ into (\ref{f2.8}), we obtain
\begin{equation}\label{f2.11}
\begin{aligned}
{{\bm{x}}}(h+&1)\\
=&(\bm{A}+\bm{b}_u\bm{K})^{2}{\bm{x}}(h-1)+\bm{b}_dd(h)\\
&+\bm{P}(\begin{bmatrix}-\bm{I}_1\bm{K}_p\bm{b}_d d(h)\cr -\bm{I}_2\bm{K}_p\bm{b}_d d(h)\end{bmatrix}+\begin{bmatrix}-\bm{I}_2\bm{K}_p\bm{b}_d d(h+1) \\ -\bm{I}_1\bm{K}_p\bm{b}_d d(h-1) \end{bmatrix})\\
=&(\bm{A}+\bm{b}_u\bm{K})^{2}{\bm{x}}(h-1)-\bm{P}\bm{K}_p\bm{b}_d d(h)\\
&-\bm{b}_u\bm{I}_2\bm{K}_p\bm{b}_d d(h+1)+\bm{b}_dd(h)\\
=&(\bm{A}+\bm{b}_u\bm{K})^{2}{\bm{x}}(h-1)-\bm{P}\bm{P}^{-1}\bm{b}_d d(h)\\
&+\bm{b}_dd(h)-\bm{b}_u\bm{I}_2\bm{K}_p\bm{b}_d d(h+1)\\
=&(\bm{A}+\bm{b}_u\bm{K})^{2}{\bm{x}}(h-1)-\bm{b}_u\bm{I}_2\bm{K}_p\bm{b}_d d(h+1).
\end{aligned}
\end{equation}
By plugging the above equation into the regulated output equation in the system (\ref{f2.1}) and using ${\bm{c}_o}{\bm{b}_u}=0$, the regulated output ${y}_{o}(h+1)$ is derived as
\begin{equation}\label{f3.8}
\begin{aligned}
{y}_{o}&(h+1)=\bm{c}_o{\bm{x}}(h+1)\\
&=\bm{c}_o((\bm{A}+\bm{b}_u\bm{K})^{2}{\bm{x}}(h-1)-\bm{b}_u\bm{I}_2\bm{K}_p\bm{b}_d d(h+1))\\\
&=\bm{c}_o(\bm{A}+\bm{b}_u\bm{K})^{2}{\bm{x}}(h-1)-\bm{c}_o\bm{b}_u\bm{I}_2\bm{K}_p\bm{b}_d d(h+1)\\
&=\bm{c}_o(\bm{A}+\bm{b}_u\bm{K})^{2}{\bm{x}}(h-1),
\end{aligned}
\end{equation}
which indicates that the disturbance in $y_o(h+1)$ are removed.

In the case of $j>h+1$, along the similar derivation with (\ref{f2.11}), one can have
\begin{equation}\label{f4.10}
\begin{aligned}
{{\bm{x}}}(j)=(\bm{A}+\bm{b}_u\bm{K})^{j-h+1}{\bm{x}}&(h-1)\notag\\
&-\bm{b}_u\bm{I}_2\bm{K}_p\bm{b}_d d(j-1).
\end{aligned}
\end{equation}
Similar to (\ref{f3.8}), there holds
\begin{align}\label{f4.11}
{y}_{o}(j)=\bm{c}_o(\bm{A}+\bm{b}_u\bm{K})^{j-h+1}{\bm{x}}(h-1).
\end{align}
It is evident that the disturbance in $y_o(j)$ is eliminated either.
\end{IEEEproof}

\begin{figure}[h!]\centering
	\includegraphics[width=8.5cm]{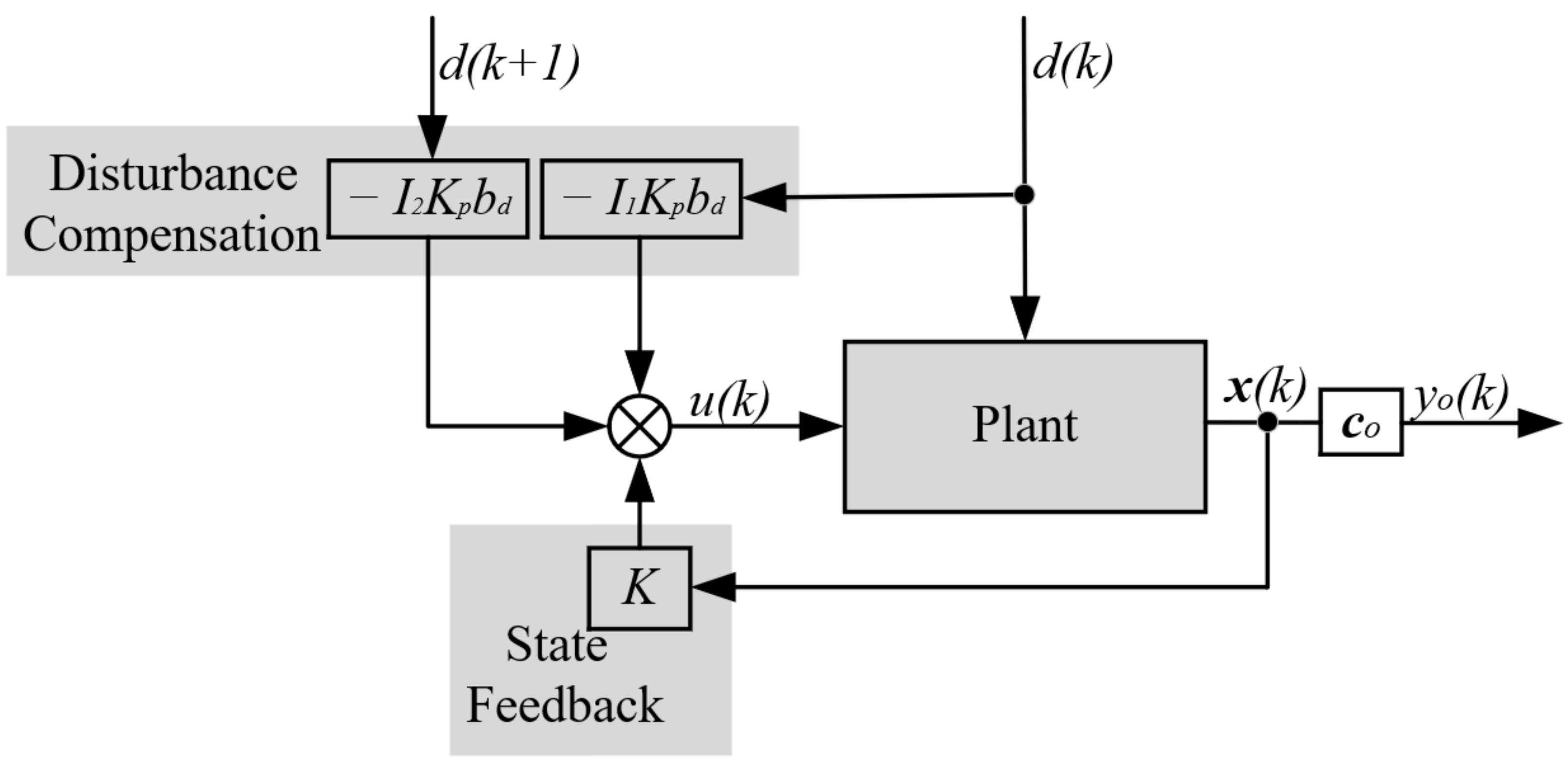}
	\caption{Configuration of the proposed mismatched known disturbance rejection control method.}\label{fig_1}
\end{figure}

The proof of Theorem \ref{the:1} is completed.
The configuration of the proposed mismatched known disturbance rejection control method proposed in Theorem \ref{the:1} is shown in Fig. \ref{fig_1}.

\begin{rem}\label{rem.3}
It is worth stressing that the controller proposed in Theorem \ref{the:1} can eliminate the disturbance in the regulated output immediately and completely, rather than in the steady state. The result is new to the best of our knowledge.
\end{rem}

\begin{rem}\label{rem.4}
Note that the disturbance compensation part $u_{d}(k)$ in (\ref{max.1}) is also suitable for matched disturbance. For the matching case (i.e., $b_u= \lambda b_d$, $ \lambda \in R$), the disturbance compensation part (\ref{3.4}) degenerates to $u_{d}(k)=-1/\lambda d(k)$, which is consistent with standard ESOBC control law \cite{5572931}.
\end{rem}

\section{Unknown Disturbance Rejection Controller}\label{sec4}
In this section, we will deal with the mismatched disturbance rejection for the system (\ref{f2.1}), where the disturbance is unknown. Similar to Assumption 2 in \cite{yang2012nonlinear} and \cite{2012Generalized}, ${d}(k)$ in \eqref{f2.1} satisfies the following assumption.
\begin{assumption}\label{assumption:3}
$d(k)$ is bounded and $\lim_{k \to \infty}{d}(k) = {c}$.
\end{assumption}

Denote
\begin{align}\label{2.2}
{h}(k)={d}(k+1)-{d}(k).
\end{align}
From Assumption 2, it is evident that $h(k)$ is bounded and $\lim_{k \to \infty}{h}(k)={0}$.

Due to the disturbance is unknown, it is necessary to design a generalized ESO to estimate and then try to reject the unknown disturbance.
\subsection{Generalized ESO Design and Analysis}

Introducing an extended variable \cite{5572931}
\begin{eqnarray}
\begin{array}{ll}
x_{3}(k)=d(k)\notag\\
\end{array}
\end{eqnarray}
to the system (\ref{f2.1}) and combine (\ref{2.2}) yield the extended system equation
\begin{equation}\label{321q.1}
\left\{
\begin{aligned}
&\bar{\bm{x}}(k+1)={\bar{\bm{A}}\bar{\bm{x}}}(k)+\bar{\bm{b}}_u u(k) +\bm{E} h(k),\\
&\bm{y}(k)=\bar{\bm{C}}\bm{x}(k)
\end{aligned}
\right.
\end{equation}

where
\begin{align}
&\bar{\bm{x}}(k)=\begin{bmatrix}\bm{x}(k)\cr x_{3}(k)\end{bmatrix},\notag\\
&\bar{\bm{A}}=\begin{bmatrix}\bm{A}_{2\times 2}&(\bm{b}_d)_{2\times 1}\cr 0&1\end{bmatrix}_{3\times3},\notag\\
&\bar{\bm{b}}_u=\begin{bmatrix}(\bm{b}_u)_{2\times1}\cr 0\end{bmatrix}_{3\times1},\notag\\
&{\bm{E}}=\begin{bmatrix}(\bm{0})_{2\times1}\cr 1\end{bmatrix}_{3\times1},\notag\\
&\bar{\bm{C}}=\begin{bmatrix}(\bm{C}_m)_{r\times 2}&\bm{0}_{r\times1}\end{bmatrix}_{r\times3}.\notag
\end{align}

For system (\ref{321q.1}), the ESO is designed as follows:
\begin{equation}\label{ESO}
\left\{
\begin{aligned}
&{\hat{\bar{\bm{x}}}(k+1)={\bar{\bm{A}}\hat{\bar{\bm{x}}}}(k)+\bar{\bm{b}}_u u(k)}+\bar{\bm{L}}(\bm{y}(k)-\hat{\bar{\bm{y}}}(k)),\\
&\hat{\bar{\bm{y}}}(k)={\bar{\bm{C}} \hat{\bar{\bm{x}}}(k)}
\end{aligned}
\right.
\end{equation}
where $\hat{\bar{\bm{x}}}(k)=\begin{bmatrix}{\hat{\bm{x}}(k)}^T,$ $\hat{x}_{3}(k)\end{bmatrix}^T$, $\hat{\bm{x}}(k)$, and $\hat{x}_{3}(k)$ are the estimates of $\bar{\bm{x}}(k)$, $\bm{x}(k)$, and ${x}_{3}(k)$  in (\ref{321q.1}), respectively. Matrix $\bar{\bm{L}}$  is the observer gain remains to be determined.

\subsection{Control Law Design}
Due to the disturbance is unknown, the control law $(\ref{max.1})$ is invalid for disturbance rejection. To this end, we will modify it by replacing $u_d(k)$ in Theorem \ref{the:1} with the below
\begin{align}\label{k_d00}
{u_{d}(k)}=-\bm{I}_1\bm{K}_p\bm{b}_d \hat d(k)-\bm{I}_2\bm{K}_p\bm{b}_d \hat d(k+1)
\end{align}
where $\hat{d}(k)$  is the estimation of $d(k)$.\\
Given the causal restriction of the controller, $u_d(k)$ in (\ref{k_d00}) is modified as
\begin{equation}\label{k_d2}
\begin{aligned}
{u_{d}(k)}&=-\bm{I}_1\bm{K}_p\bm{b}_d \hat d(k)-\bm{I}_2\bm{K}_p\bm{b}_d \hat d(k)\\
&=-(\bm{I}_1+\bm{I}_2)\bm{K}_p\bm{b}_d \hat d(k)=K_d \hat d(k).
\end{aligned}
\end{equation}
Thus, a new control law is obtained as
\begin{align}\label{max.2}
u(k)={u}_{x}(k)+{u_{d}(k)}=\bm{K}\bm{x}(k)+K_d \hat d(k),
\end{align}
which is applicable in the case that the state is available. If the state is unavailable, the control law (\ref{max.2}) can be modified as
\begin{align}\label{max.3}
u(k)=\bm{K}\hat{\bm{x}}(k)+K_d \hat d(k).
\end{align}

Until now, we have no idea whether the proposed controller (\ref{max.2}) and (\ref{max.3}) can stabilize the system (\ref{f2.1}), or/and eliminate the disturbance in regulated output or not. Hence, what follows is stability and disturbance rejection analysis.

\subsection{Stability and Disturbance Rejection Analysis}
Before analysis, some useful lemmas are prepared as follows.

\begin{lemma}\label{lemma.2}
The following linear system
\begin{align}\label{f0.1}
{x}(k+1)= {A x}(k) + {B}u(k)
\end{align}
is asymptotically stable if ${A}$ is Schur and u(k) is bounded as well as $\lim_{k \to \infty}u(k)=0$.
\end{lemma}
\begin{IEEEproof}
Let $X(z)$ and $U(z)$ are z-transform of $x(k)$ and $u(k)$, respectively. According to the final value theorem \cite{lathi2005linear}, if $(z-1)X(z)$ has no poles outside the unit circle, then there holds
\begin{align}\label{A.1}
\lim_{k \to \infty}{x}(k)=\lim_{z \to 1}(z-1)X(z). 
\end{align}
The  z-transform  of (\ref{f0.1}) is 
\begin{align}\label{100.1}
X(z)=(z{I}-{A})^{-1}BU(z).
\end{align}
Further, we have 
\begin{align}\label{100.2}
(z-1)X(z)=(z-1)(z{I}-{A})^{-1}BU(z).
\end{align}
It follows from the boundedness of $u(k)$ and $\lim_{k \to \infty}{u}(k)={0}$ that
 $(z-1)U(z)$ has no poles outside the unit circle. Applying the final value theorem again yields \begin{align}
\lim_{z \to 1}(z-1)U(z)=\lim_{k \to \infty}u(k)=0.  \label{uk0}\end{align}
Since $A$ is Schur, all poles of $(z{I}-{A})^{-1}$ are inside the unit circle and 
$\lim_{z\to 1} (zI-A)^{-1}$ exists, which together with \eqref{uk0} shows that
\begin{align}
\lim_{z \to 1}(z{I}-{A})^{-1}B(z-1)U(z)=0. \label{twolim}
\end{align}
In view of \eqref{100.2} and \eqref{twolim}, 
 \begin{align} \lim_{z \to 1}(z-1)X(z)=0. \end{align}
Now from \eqref{A.1}, it is immediate to obtain that 
\begin{align}\label{0.1}
\lim_{k \to \infty}{x}(k)={0}.
\end{align}
\end{IEEEproof}

\begin{lemma}\label{lemma.3}
For system (\ref{f0.1}), $\lim_{k \to \infty}x(k)={({I}-{A})}^{-1}{B}U_c$ if $A$ is Schur and $\lim_{k \to \infty}u(k)=U_c \neq 0$.
\end{lemma}
\begin{IEEEproof}Along the similar reasoning with Lemma \ref{lemma.2}, Lemma \ref{lemma.3} can be derived. The proof is thus omitted.
\end{IEEEproof}

Let
\begin{align}\label{5.4}
\bm{e}(k)=\begin{bmatrix}\bm{e}_x(k)\cr e_d(k)\end{bmatrix}= \begin{bmatrix}\hat{\bm{x}}(k)-\bm{x}(k)\cr \hat{d}(k)-{d}(k)\end{bmatrix}.
\end{align}

By using (\ref{321q.1}), (\ref{ESO}) and (\ref{5.4}), we obtain the following estimation error equation
\begin{align}\label{5.6}
\bm{e}(k+1)= (\bar{\bm{A}}-\bar{\bm{L}}\bar{\bm{C}})\bm{e}(k)-\bm{E}h(k).
\end{align}

Under Assumption $\ref{assumption:3}$, if selects the observer gain matrix $\bar{\bm{L}}$ in (\ref{5.6}) such that $\bar{\bm{A}}-\bar{\bm{L}}\bar{\bm{C}}$ is Schur matrix, then according to Lemma \ref{lemma.2} the observer error system (\ref{5.6}) is asymptotically stable and
\begin{align}\label{5.12}
\lim_{k \to \infty}\bm{e}(k)=\bm{0}
\end{align}
for bounded $h(k)$ and $\lim_{k \to \infty}{h}(k)={0}$.

\subsubsection{In the Case of Known State}
If the state is accessible, the control law can be given as (\ref{max.2}). We will analyze the performance of the controller in the below.

\begin{theorem}\label{theorem.2}
Under Assumption \ref{assumption.1} and \ref{assumption:3}, $(\bm{A}, \bm{b}_u)$ is controllable, and $(\bar{\bm{A}}, \bar{\bm{C}})$ is observable.

1) Select the observer gain matrix $\bar{\bm{L}}$ in (\ref{ESO}) such that $\bar{\bm{A}}-\bar{\bm{L}}\bar{\bm{C}}$ is Schur matrix \cite{lathi2005linear};

2) Choose the state feedback gain matrix $\bm{K}$ in (\ref{max.1}) such that $\bm{A}+\bm{b}_u\bm{K}$ is Schur matrix;

Under the composite control law (\ref{max.2}), there hold
\begin{itemize}
\item The closed-loop system (\ref{f2.1}) is asymptotically stable;
\item The disturbance ${d}(k)$ can be eliminated from the steady-state regulated output ${y}_{o}(k)$.
\end{itemize}

\end{theorem}

\begin{IEEEproof}
In view of (\ref{f2.1}), (\ref{max.2}) and (\ref{5.6}), the system (\ref{f2.1}) under the control law (\ref{max.2}) evolves as
\begin{equation}\label{com.1}
\begin{aligned}
\bm{x}(k+1)=(\bm{A}+\bm{b}_u\bm{K})\bm{x}(k)+&\bm{b}_u\bar{\bm{K}}\bm{e}(k)\\
&+(\bm{b}_d+\bm{b}_uK_d) d(k)
\end{aligned}
\end{equation}
where $\bar{\bm{K}}=\begin{bmatrix}\bm{0}_{1 \times 2}&K_d\end{bmatrix}$.

Since $\bm{A}+\bm{b}_u\bm{K}$ is Schur matrix. Thus, it can be concluded from Lemma \ref{lemma.3} and  Assumption $\ref{assumption:3}$ that the state $\bm{x}(k)$ in the closed-loop system (\ref{com.1}) is convergent.
It can be derived from (\ref{com.1}) that
\begin{equation}\label{5.10}
\begin{aligned}
\lim_{k \to \infty}\bm{x}(k)=\lim_{k \to \infty}&{(\bm{I}-(\bm{A}+\bm{b}_u\bm{K}))}^{-1}[\bm{b}_u\bar{\bm{K}}\bm{e}(k)\\
&+(\bm{b}_d+\bm{b}_uK_d)d(k)].
\end{aligned}
\end{equation}
Combining (\ref{f2.1}) with (\ref{5.10}) gives
\begin{equation}\label{5.11}
\begin{aligned}
\lim_{k \to \infty}y_o(k)=\lim_{k \to \infty}\bm{c}_o&\bm{x}(k)\\
=\lim_{k \to \infty}\bm{c}_o&{(\bm{I}-(\bm{A}+\bm{b}_u\bm{K}))}^{-1}[\bm{b}_u\bar{\bm{K}}\bm{e}(k)\\
&+(\bm{b}_d+\bm{b}_uK_d)d(k)].
\end{aligned}
\end{equation}
Substituting (\ref{5.12}) into above equation yields
\begin{equation}\label{5.31}
\begin{aligned}
\lim_{k \to \infty}y_o(k)\lim_{k \to \infty}\bm{c}_o{(\bm{I}-(\bm{A}+\bm{b}_u\bm{K}))}^{-1}(\bm{b}_d+\bm{b}_uK_d)d(k).
\end{aligned}
\end{equation}

Using (\ref{k_d2}) and (\ref{5.31}) yields
\begin{equation}\label{4.21}
\begin{aligned}
\lim_{k \to \infty}y_o(k)=&\lim_{k \to \infty}\bm{c}_o{(\bm{I}-(\bm{A}+\bm{b}_u\bm{K}))}^{-1}\\
&\times(\bm{b}_d-\bm{b}_u(\bm{I}_1+\bm{I}_2)\bm{K}_p\bm{b}_d)d(k)\\
=&\lim_{k \to \infty}\bm{c}_o{(\bm{I}-(\bm{A}+\bm{b}_u\bm{K}))}^{-1}\\
&\times(\bm{P}\bm{P}^{-1}-\bm{b}_u(\bm{I}_1+\bm{I}_2)\bm{P}^{-1})\bm{b}_d d(k)\\
=&\lim_{k \to \infty}\bm{c}_o{(\bm{I}-(\bm{A}+\bm{b}_u\bm{K}))}^{-1}\\
&\times\begin{bmatrix}0,(\bm{A}+\bm{b}_u\bm{K}-\bm{I}){\bm{b}_u}\end{bmatrix}\bm{P}^{-1}\bm{b}_d d(k)\\
=&\lim_{k \to \infty}-\bm{c}_o\bm{b}_u\bm{I}_2\bm{P}^{-1}\bm{b}_d d(k).
\end{aligned}
\end{equation}
Substituting the ${\bm{c}_o}{\bm{b}_u}=0$ into above equation yields
\begin{align}\label{5.13}
\lim_{k \to \infty}y_o(k)=0
\end{align}
which means that the disturbance $d(k)$ can be removed from the steady-state regulated output $y_o(k)$ under the composite control law (\ref{max.2}).

\end{IEEEproof}

\begin{figure}[!t]\centering
	\includegraphics[width=8.5cm]{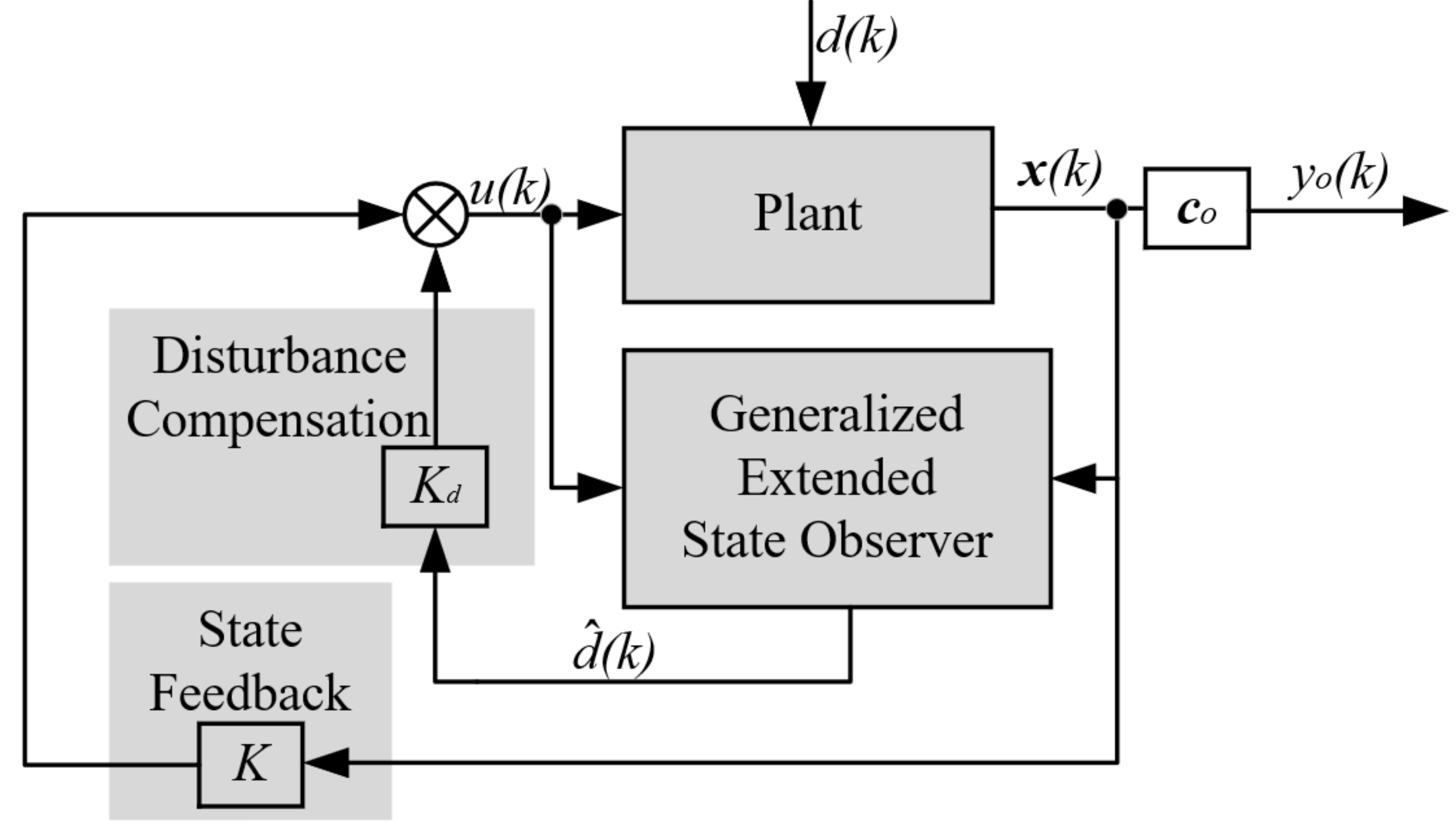}
	\caption{Configuration of the mismatched unknown disturbance rejection control method.}\label{fig_2}
\end{figure}

The proof of Theorem \ref{theorem.2} is completed.
The configuration of the proposed mismatched unknown disturbance rejection control method for the state available is shown in Fig. \ref{fig_2}.

\subsubsection{In the Case of Unknown State}
If the state is unavailable, the control law can be given as (\ref{max.3}). We will analyze the performance of the controller in the below.
\begin{theorem}\label{theorem.3}
Under Assumption \ref{assumption.1} and \ref{assumption:3}, $(\bm{A}, \bm{b}_u)$ is controllable, and $(\bar{\bm{A}}, \bar{\bm{C}})$ is observable.

1) Select the observer gain matrix $\bar{\bm{L}}$ in (\ref{ESO}) such that $\bar{\bm{A}}-\bar{\bm{L}}\bar{\bm{C}}$ is Schur matrix;

2) Choose the state feedback gain matrix $\bm{K}$ in (\ref{max.1}) such that $\bm{A}+\bm{b}_u\bm{K}$ is Schur matrix;

Under the composite control law (\ref{max.3}), there hold
\begin{itemize}
\item The closed-loop system (\ref{f2.1}) is asymptotically stable;
\item The disturbance ${d}(k)$ can be eliminated from the steady-state regulated output ${y}_{o}(k)$.
\end{itemize}

\end{theorem}
\begin{IEEEproof}
In view of (\ref{f2.1}), (\ref{max.3}) and (\ref{5.6}), the system (\ref{f2.1}) under the control law (\ref{max.3}) evolves as
\begin{equation}\label{com.2}
\begin{aligned}
\bm{x}(k+1)=(\bm{A}+\bm{b}_u\bm{K})\bm{x}(k)+&\bm{b}_u\tilde{\bm{K}}\bm{e}(k)\\
&+(\bm{b}_d+\bm{b}_uK_d) d(k)
\end{aligned}
\end{equation}
where $\tilde{\bm{K}} = [\bm{K},K_d]$.

Since $\bm{A}+\bm{b}_u\bm{K}$ are Schur matrix. According to Lemma \ref{lemma.3}, the state $\bm{x}(k)$ in close-loop system (\ref{com.2}) under Assumption $\ref{assumption:3}$ is convergent. From (\ref{com.2}), it can be derived that
\begin{equation}\label{5.16}
\begin{aligned}
\lim_{k \to \infty}\bm{x}(k)=\lim_{k \to \infty}&{(\bm{I}-(\bm{A}+\bm{b}_u\bm{K}))}^{-1}\\
&\times[\bm{b}_u\tilde{\bm{K}}\bm{e}(k)+(\bm{b}_d+\bm{b}_uK_d)d(k)].
\end{aligned}
\end{equation}
It follows from (\ref{f2.1}), (\ref{5.12}), and (\ref{5.16}) that the regulated output can be represented as
\begin{equation}\label{5.17}
\begin{aligned}
\lim_{k \to \infty}&y_o(k)=\lim_{k \to \infty}\bm{c}_o\bm{x}(k)\\
=&\lim_{k \to \infty}\bm{c}_o{(\bm{I}-(\bm{A}+\bm{b}_u\bm{K}))}^{-1}(\bm{b}_d+\bm{b}_uK_d)d(k).
\end{aligned}
\end{equation}
With the same lines as (\ref{4.21}), the same result as (\ref{5.13}) can be obtained from (\ref{5.17}).
\end{IEEEproof}

The proof of Theorem \ref{theorem.3} is completed.

\begin{rem}
From Remark \ref{rem.3}, the residual disturbance in the regulated output under the modified controller is essentially induced by the error of observer. Consequently, the modified controller can eliminate the disturbance in the regulated output
completely if the observer can provide the exact disturbance value at every instant.
\end{rem}

\section{Simulations}\label{sec5}
\subsection{Example with Known Disturbance}
In the case of known disturbance, the designed control law in this paper is compared with the GESOBC\cite{2012Generalized} to illustrate the effect of disturbance compensation.

Consider the system (\ref{f2.1}) with the following parameters
\begin{equation}
\begin{aligned}
&\bm{A}=\begin{bmatrix}1&0.01\cr -0.02&0.99\end{bmatrix}, \quad \bm{b}_u=\begin{bmatrix}0\cr 0.01\end{bmatrix},\quad \bm{b}_d=\begin{bmatrix}0.01\cr 0\end{bmatrix},\notag\\
&\bm{C}_m=\begin{bmatrix}1&0\cr 0&1\end{bmatrix} , \quad \bm{c}_o=\begin{bmatrix}1&0\end{bmatrix}.
\end{aligned}
\end{equation}

\begin{figure}[h!]\centering
	\includegraphics[width=8.5cm,height=7cm]{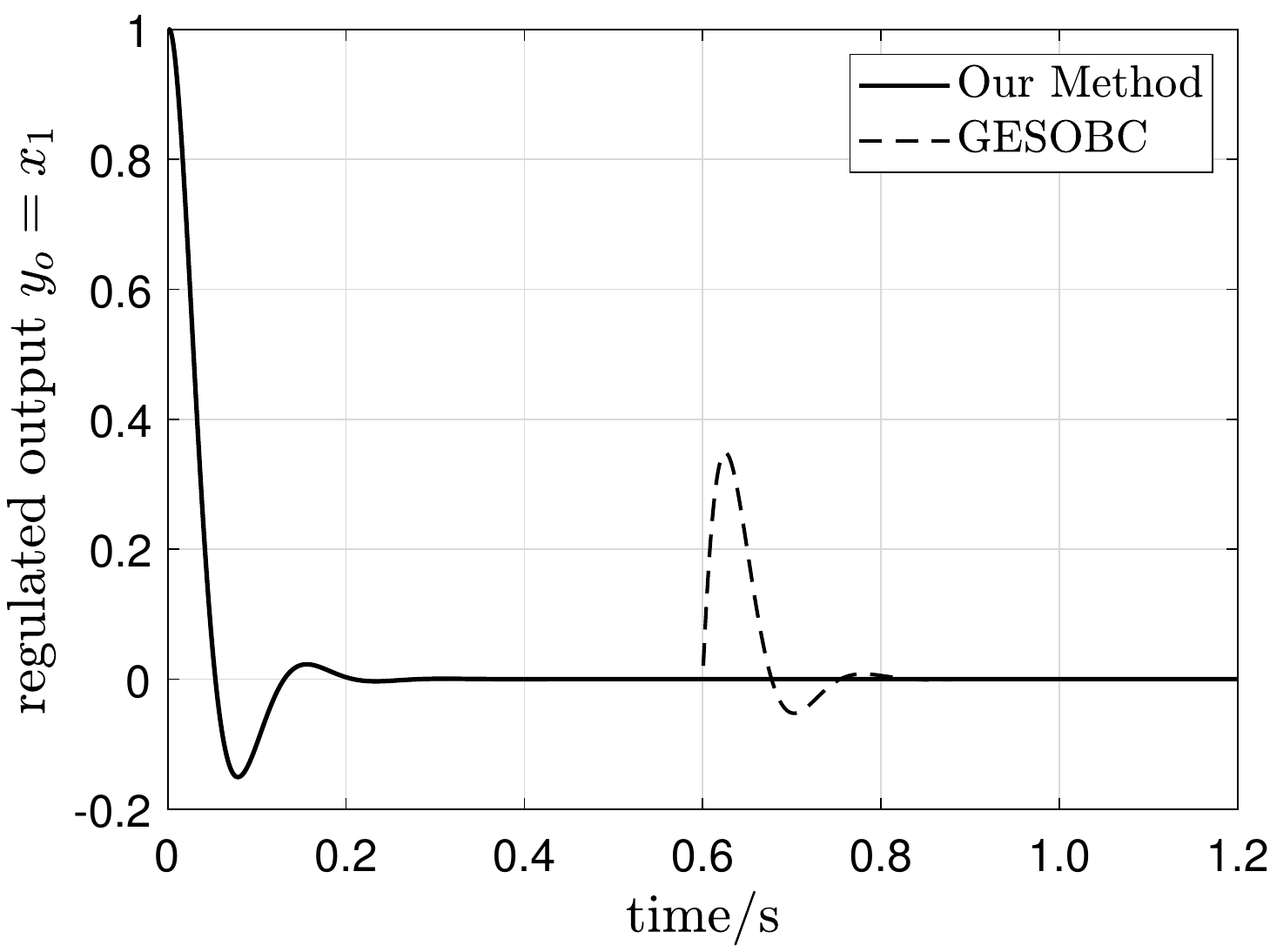}
	\caption{Simulation trajectories of the regulated output $y_o=x_1$ generated by two methods.}\label{fig_3}
\end{figure}

\begin{figure}[!t]\centering
	\includegraphics[width=8.5cm,height=7cm]{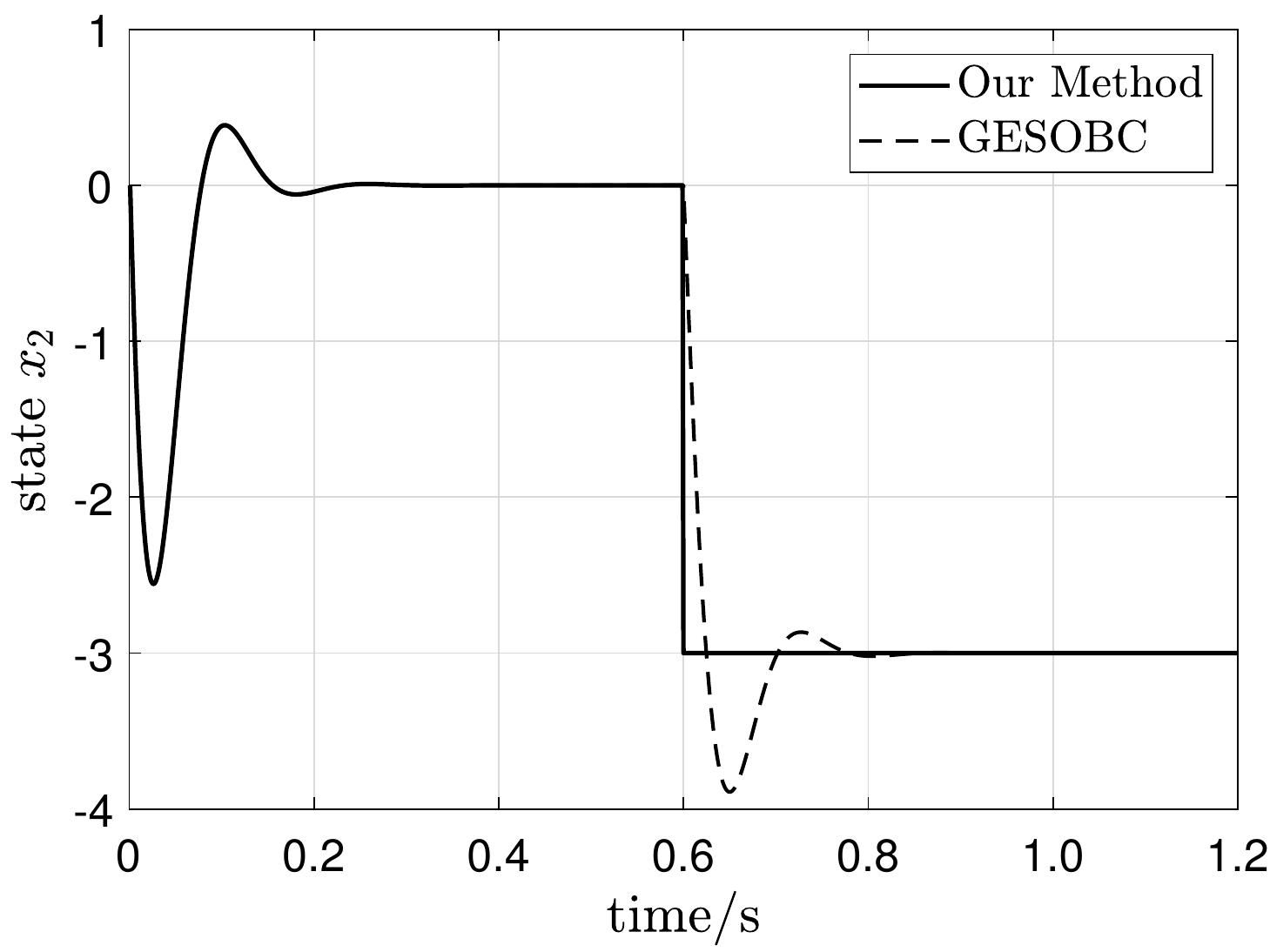}
	\caption{Simulation trajectories of state $x_2$ generated by two methods.}\label{fig_4}
\end{figure}

The state feedback gain matrix in (\ref{max.1}) is chosen as $\bm{K} =\begin{bmatrix}-20&-4\end{bmatrix}$ such that the poles of closed-loop system regardless of the disturbance are $0.9750 + 0.0397i$ and $0.9750 - 0.0397i$. Two coefficients of the disturbance compensation part $u_d(k)$ can be calculated according to (\ref{max.1}), giving as $-\bm{I}_1\bm{K}_p\bm{b}_d =95$ and$-\bm{I}_2\bm{K}_p\bm{b}_d = -100$. According to \cite{2012Generalized}, the disturbance compensation gain of the GESOBC method is calculated as $K_d=-5$. The initial state of the system  is $x(0) =\begin{bmatrix}1&0\end{bmatrix}'$. The known disturbance $d=3$ acts on the system from $t = 0.6$ $s$. The controller aims to remove the disturbance from the regulated output $y_o$. The simulation trajectories of the regulated output $y_o=x_1$ and the state $x_2$ are shown in Figs. \ref{fig_3}-\ref{fig_4}.

In Fig. \ref{fig_3}, it can be clearly seen that the proposed method is more effective in disturbance rejection than GESOBC because the proposed method has completely eliminated the disturbance in the regulated output $y_o=x_1$. The state $x_2$ (shown in Fig. \ref{fig_4}) reaches the steady state without any oscillation induced by the disturbance.

\subsection{Example with Unknown Disturbance}
Consider the permanent-magnet direct current (PMDC) motor system \cite{2002Analysis} is described by 
\begin{align}
V_a(t)=&R_ai_a(t)-L_a\frac{di_a(t)}{dt}+E_b(t),\notag\\
&E_b(t)=K_bw_m(t),\notag\\
&T_m(t)=K_ti_a(t),\label{model1}\\
&K_t=K_b=k,\notag\\
T_m(t)=&J_m\frac{dw_m(t)}{dt}+B_mw_m+T_L(t),\notag
\end{align}
where $V_a(t)$ is the armature voltage $(V)$, $d(t)$ = $T_L(t)$ is load torque $(Nm)$, $w_m$  is angular speed $(rad/s)$; $i_a(t)$  is armature current$(A)$; $R_a$, $L_a$ respectively the armature resistance $(\Omega)$ and armature inductance $(H)$; $T_m$ is motor torque $(Nm)$; $E_b$ is back emf $(V)$; $K_b$ back emf constant $(Vs/rad)$; $K_t$ is torque constant $(Nm/A)$; $J_m$ is rotor inertia $(kgm^2)$; $B_m$ is  the viscous friction coefficient $(Nms/rad)$; $\theta$ is angular position of rotor shaft$(rad)$.
The parameters in the system (\ref{model1}) are listed in Table \ref{table_1} \cite{2015Effective}.

Let $\bm x(t)=[ w_m\quad i_a]'$, $\bm{y}(t)=[ w_m\quad i_a]'$, ${y}_o(t)=w_m$, $u(t) =V_a(t)$, so the state-space equation of the PMDC motor model as
\begin{numcases}{}
{\dot{\bm{{x}}}(t)=\bar{\bm{A}}\bm{x}(t)+\bar{\bm{b}}_u{u}(t)}+\bar{\bm{b}}_d{d}(t),\notag\\
\bm{y}(t)=\bm{Cx}(t),\label{321p.1}\\
{y}_o(t)=\bm{c}_o\bm{x}(t),\notag
\end{numcases}
where
\begin{equation}\label{321v.1}
\begin{aligned}
&{\bm{A}}=\begin{bmatrix}-\dfrac{B_m}{J_m}&\dfrac{k}{J_m}\\[2.5mm]-\dfrac{k}{L_{a}}&-\dfrac{R_a}{L_{a}}\end{bmatrix}, \bar{\bm{b}}_u=\begin{bmatrix}0\\[2.5mm] \dfrac{1}{L_{a}}\end{bmatrix}, &\bar{\bm{b}}_d=\begin{bmatrix}-\dfrac{1}{J_m}\\[2.5mm] 0 \end{bmatrix},\notag\\
&\bm{C}=\begin{bmatrix}1&0\cr 0&1\end{bmatrix}, \quad \bm{c}_o=\begin{bmatrix}1&0\end{bmatrix}.
\end{aligned}
\end{equation}

\begin{table}
\caption{Parameters of PMDC Motor}
\label{table}
\setlength{\tabcolsep}{3pt}
\begin{tabular}{|p{25pt}|p{75pt}|p{115pt}|}
\hline
Symbol& 
Explanation& 
Value \\
\hline
  $P$   & Shaft Power      &  $5kW$ \\
  $V_a$ & Armature Voltage  & $240V$ \\
  $R_a$ &Armature Resistance &$0.5\Omega$\\
  $L_a$ &Armature Inductance &$0.012H$\\
  $J_m$ &Total Inertia &$0.00471kgm^2$\\
  $B_m$ &Viscous Friction Coef. &$0.002Nms/rad$\\
  $K_t$ &Torque Constant &$0.5Nm/A$\\ 
  $K_b$ &Back Emf Constant &$0.5Vs/rad$\\
\hline
\end{tabular}
\label{table_1}
\end{table}

The disturbance in PMDC motor system is the load torque. The control strategy in this paper can be utilized to study the speed control on PMDC motor under variable loads. It can be noted from (\ref{321p.1}) and (\ref{321v.1}) that the disturbance in the PMDC motor system is mismatched. The control objective is to eliminate the disturbance in the output angular velocity under variable loads.

Under the method proposed in this paper, the sampling period is selected as $T=0.001s$ and disturbance change is $5Nm$ loaded from $0.6s$ action. The state feedback gain matrix designed as $\bm{K}=\begin{bmatrix}-0.5&-4\end{bmatrix}$. It can be calculated from (\ref{k_d2}) that $K_d=5.3$.
The the observer gain of the generalized ESO is chosen as
\begin{align}\label{123}
\bar{\bm{L}}=\begin{bmatrix}0.3&0.1\\0.1&0.8\\-0.2&-0.05\end{bmatrix}.
\end{align}
\par
\begin{figure}[h]\centering
	\includegraphics[width=8.5cm,height=7cm]{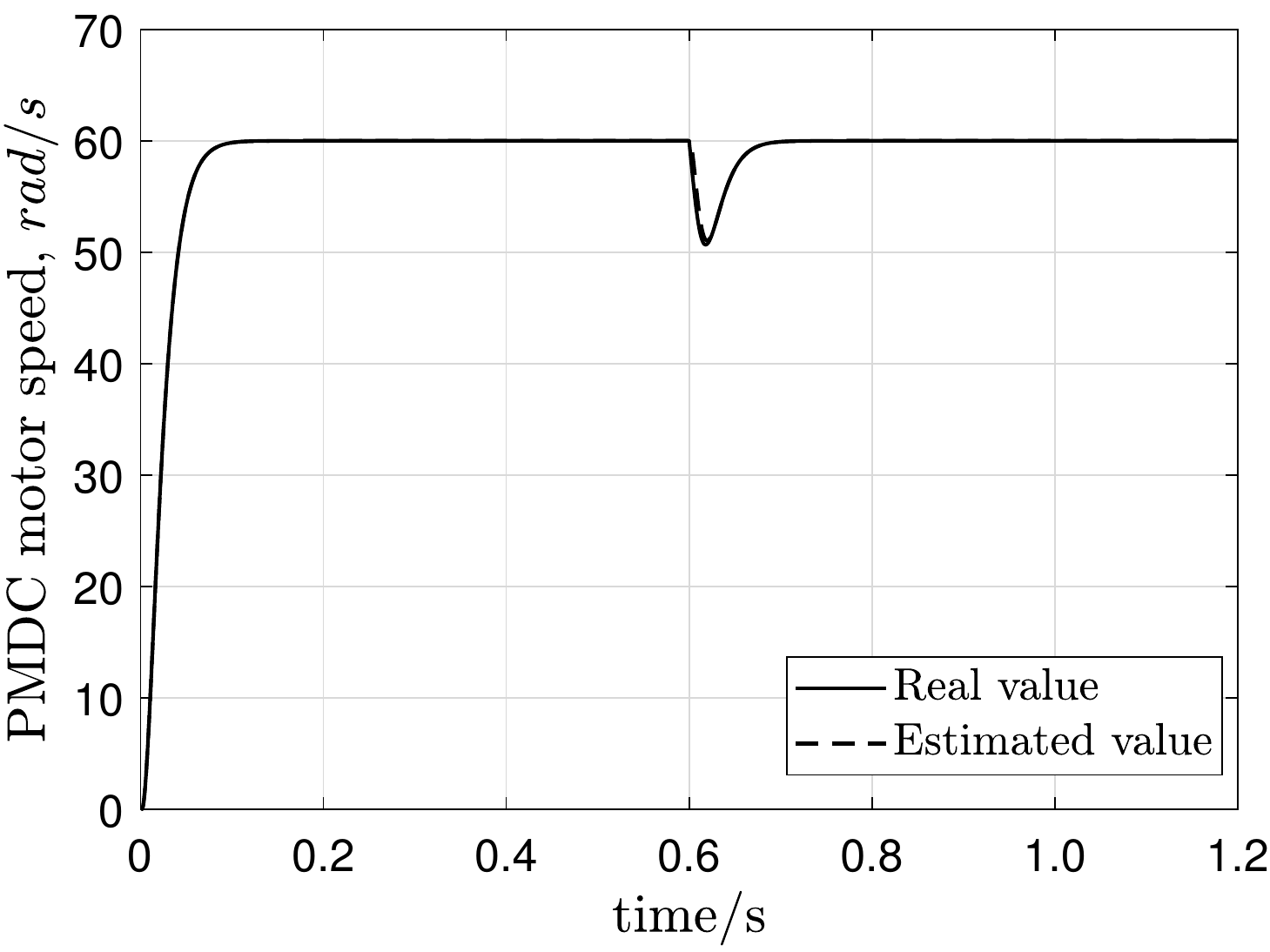}
	\caption{Simulation trajectories of the PMDC motor speed and its estimated value.}\label{fig_5}
\end{figure}

\begin{figure}[h]\centering
	\includegraphics[width=8.5cm,height=7cm]{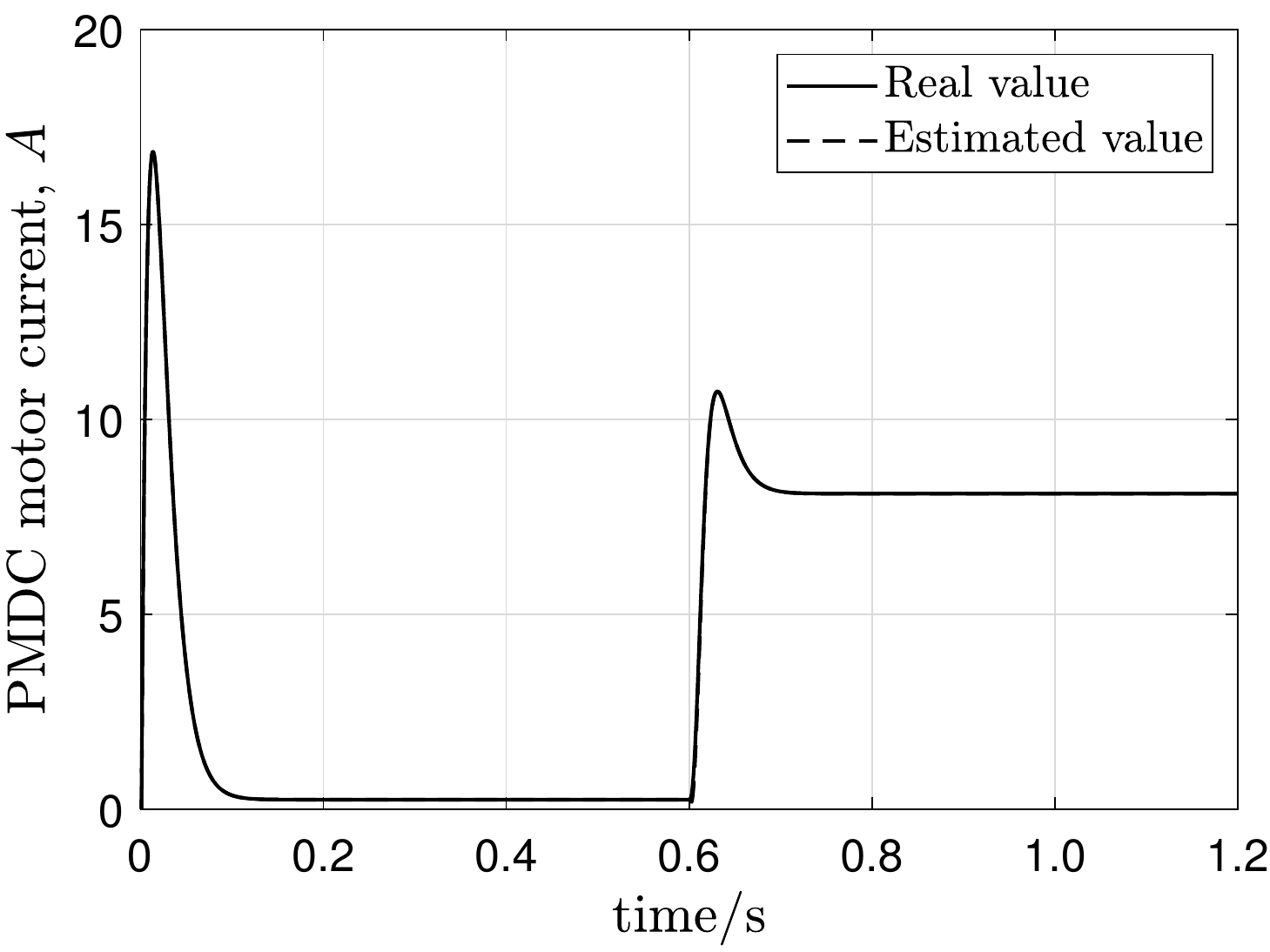}
	\caption{Simulation trajectories of the PMDC motor current and its estimated value.}\label{fig_6}
\end{figure}

\begin{figure}[h]\centering
	\includegraphics[width=8.5cm,height=7cm]{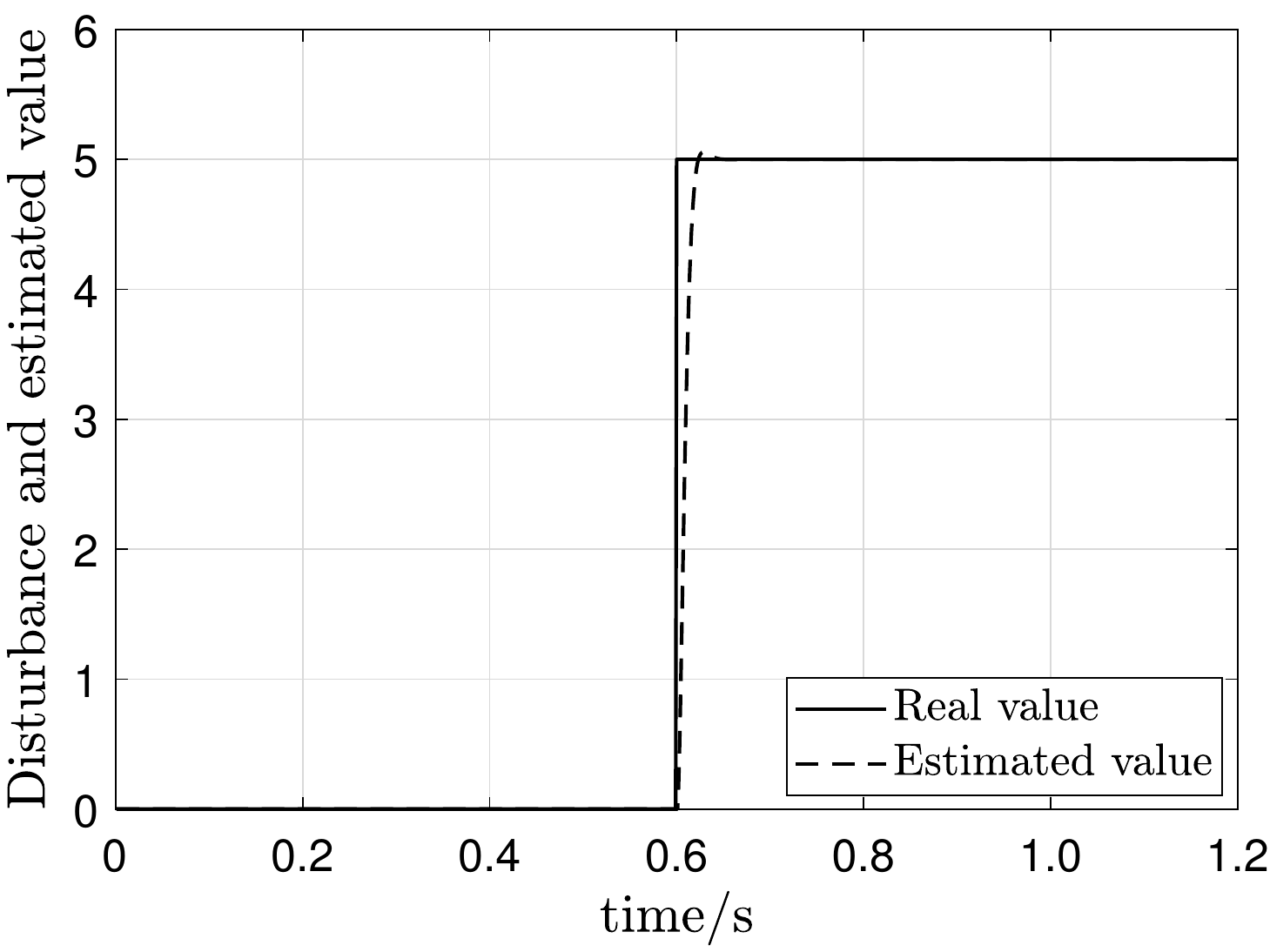}
	\caption{Trajectories of the disturbance and its estimated value.}\label{fig_7}
\end{figure}

 The simulation trajectories of the PMDC motor speed output real and estimated value under our approach are shown in Fig. \ref{fig_5}. The simulation trajectories of the real and estimated value of PMDC motor current are shown in Fig. \ref{fig_6}.  The trajectories of the real and estimated value of disturbance are shown in Fig. \ref{fig_7}.

As we can see in Fig. \ref{fig_5}-\ref{fig_6} that the maximum speed change does not exceed 10 $rad/s$, the PMDC motor current change is less than $3A$, the settling time of them is about $0.1s$, and the steady-state error disappears. Both the speed and the current vary approach to the setpoints quickly. Disturbances can be estimated accurately and timely by the generalized ESO  (see Fig. \ref{fig_7}). The results illustrate that the proposed approach has achieved expected performance in eliminating the unknown disturbance.

\section{Conclusion}\label{sec6}

In this paper, we have provided several mismatch disturbance rejection controllers for the known disturbance and the unknown disturbance. Different from the existing controllers, our disturbance compensation gain is proposed with the aid of the controllability of the system. The controller for the known disturbance can eliminate the disturbance in the regulated output immediately and completely, and the controller treating the unknown disturbance can reject the disturbance in the regulated output in the steady state. Two examples have been designed to illustrate that the proposed disturbance rejection control methods are effective.

\bibliography{ref}

\begin{thebibliography}{10}
\providecommand{\url}[1]{#1}
\csname url@samestyle\endcsname
\providecommand{\newblock}{\relax}
\providecommand{\bibinfo}[2]{#2}
\providecommand{\BIBentrySTDinterwordspacing}{\spaceskip=0pt\relax}
\providecommand{\BIBentryALTinterwordstretchfactor}{4}
\providecommand{\BIBentryALTinterwordspacing}{\spaceskip=\fontdimen2\font plus
\BIBentryALTinterwordstretchfactor\fontdimen3\font minus
  \fontdimen4\font\relax}
\providecommand{\BIBforeignlanguage}[2]{{%
\expandafter\ifx\csname l@#1\endcsname\relax
\typeout{** WARNING: IEEEtran.bst: No hyphenation pattern has been}%
\typeout{** loaded for the language `#1'. Using the pattern for}%
\typeout{** the default language instead.}%
\else
\language=\csname l@#1\endcsname
\fi
#2}}
\providecommand{\BIBdecl}{\relax}
\BIBdecl

\bibitem{1372532}
W.-H. Chen, ``Disturbance observer based control for nonlinear systems,''
  \emph{IEEE/ASME Transactions on Mechatronics}, vol.~9, no.~4, pp. 706--710,
  2004.

\bibitem{han1995extended}
J.~Han, ``The extended state observer of a class of uncertain systems,''
  \emph{Control and Decision}, vol.~10, no.~1, pp. 85--88, 1995.

\bibitem{johnson1971accomodation}
C.~Johnson, ``Accomodation of external disturbances in linear regulator and
  servomechanism problems,'' \emph{IEEE Transactions on Automatic Control},
  vol.~16, no.~6, pp. 635--644, 1971.

\bibitem{guo2014anti}
L.~Guo and S.~Cao, ``Anti-disturbance control theory for systems with multiple
  disturbances: A survey,'' \emph{ISA Transactions}, vol.~53, no.~4, pp.
  846--849, 2014.

\bibitem{5572931}
Y.~Huang, W.~Xue, and X.~Yang, ``Active disturbance rejection control:
  Methodology, theoretical analysis and applications,'' in \emph{Proceedings of
  the 29th Chinese Control Conference}, 2010, pp. 6083--6090.

\bibitem{MADONSKI201518}
R.~Mado{\'n}ski and P.~Herman, ``Survey on methods of increasing the efficiency
  of extended state disturbance observers,'' \emph{ISA Transactions}, vol.~56,
  pp. 18--27, 2015.

\bibitem{9105115}
Y.~Wang, W.~Zhang, and L.~Yu, ``A linear active disturbance rejection control
  approach to position synchronization control for networked interconnected
  motion system,'' \emph{IEEE Transactions on Control of Network Systems},
  vol.~7, no.~4, pp. 1746--1756, 2020.

\bibitem{roman2021hybrid}
R.-C. Roman, R.-E. Precup, and E.~M. Petriu, ``Hybrid data-driven fuzzy active
  disturbance rejection control for tower crane systems,'' \emph{European
  Journal of Control}, vol.~58, pp. 373--387, 2021.

\bibitem{9345466}
X.~Hu, X.~Wei, Y.~Kao, and J.~Han, ``Robust synchronization for under-actuated
  vessels based on disturbance observer,'' \emph{IEEE Transactions on
  Intelligent Transportation Systems}, pp. 1--10, 2021.

\bibitem{ding2021extended}
K.~Ding and Q.~Zhu, ``Extended dissipative anti-disturbance control for delayed
  switched singular semi-markovian jump systems with multi-disturbance via
  disturbance observer,'' \emph{Automatica}, vol. 128, p. 109556, 2021.

\bibitem{johnson1970further}
C.~Johnson, ``Further study of the linear regulator with disturbances--the case
  of vector disturbances satisfying a linear differential equation,''
  \emph{IEEE Transactions on Automatic Control}, vol.~15, no.~2, pp. 222--228,
  1970.

\bibitem{4391082}
J.~She, M.~Fang, Y.~Ohyama, H.~Hashimoto, and M.~Wu, ``Improving
  disturbance-rejection performance based on an equivalent-input-disturbance
  approach,'' \emph{IEEE Transactions on Industrial Electronics}, vol.~55,
  no.~1, pp. 380--389, 2008.

\bibitem{han2009pid}
J.~Han, ``From {PID} to active disturbance rejection control,'' \emph{IEEE
  transactions on Industrial Electronics}, vol.~56, no.~3, pp. 900--906, 2009.

\bibitem{li2014disturbance}
S.~Li, J.~Yang, W.~Chen, and X.~Chen, \emph{Disturbance observer-based control:
  methods and applications}.\hskip 1em plus 0.5em minus 0.4em\relax CRC Press,
  2014.

\bibitem{7265050}
W.-H. Chen, J.~Yang, L.~Guo, and S.~Li, ``Disturbance-observer-based control
  and related methods -- an overview,'' \emph{IEEE Transactions on Industrial
  Electronics}, vol.~63, no.~2, pp. 1083--1095, 2015.

\bibitem{2012Generalized}
S.~Li, J.~Yang, W.~H. Chen, and X.~Chen, ``Generalized extended state observer
  based control for systems with mismatched uncertainties,'' \emph{IEEE
  Transactions on Industrial Electronics}, vol.~59, no.~12, pp. 4792--4802,
  2012.

\bibitem{mohamed2007design}
Y.~A.~I. Mohamed, ``Design and implementation of a robust current-control
  scheme for a {PMSM} vector drive with a simple adaptive disturbance
  observer,'' \emph{IEEE Transactions on Industrial Electronics}, vol.~54,
  no.~4, pp. 1981--1988, 2007.

\bibitem{chwa2004compensation}
D.~Chwa, J.~Y. Choi, and J.~H. Seo, ``Compensation of actuator dynamics in
  nonlinear missile control,'' \emph{IEEE Transactions on Control Systems
  Technology}, vol.~12, no.~4, pp. 620--626, 2004.

\bibitem{chen2003nonlinear}
W.-H. Chen, ``Nonlinear disturbance observer-enhanced dynamic inversion control
  of missiles,'' \emph{Journal of Guidance, Control, and Dynamics}, vol.~26,
  no.~1, pp. 161--166, 2003.

\bibitem{isidori1985nonlinear}
A.~Isidori, \emph{Nonlinear control systems: an introduction}.\hskip 1em plus
  0.5em minus 0.4em\relax Springer, 1985.

\bibitem{castillo2018enhanced}
A.~Castillo, P.~Garc{\'\i}a, R.~Sanz, and P.~Albertos, ``Enhanced extended
  state observer-based control for systems with mismatched uncertainties and
  disturbances,'' \emph{ISA Transactions}, vol.~73, pp. 1--10, 2018.

\bibitem{ginoya2013sliding}
D.~Ginoya, P.~Shendge, and S.~Phadke, ``Sliding mode control for mismatched
  uncertain systems using an extended disturbance observer,'' \emph{IEEE
  Transactions on Industrial Electronics}, vol.~61, no.~4, pp. 1983--1992,
  2013.

\bibitem{chen2016adrc}
S.~Chen, W.~Bai, and Y.~Huang, ``{ADRC} for systems with unobservable and
  unmatched uncertainty,'' in \emph{2016 35th Chinese Control Conference
  (CCC)}.\hskip 1em plus 0.5em minus 0.4em\relax IEEE, 2016, pp. 337--342.

\bibitem{yang2012nonlinear}
J.~Yang, S.~Li, and W.-H. Chen, ``Nonlinear disturbance observer-based control
  for multi-input multi-output nonlinear systems subject to mismatching
  condition,'' \emph{International Journal of Control}, vol.~85, no.~8, pp.
  1071--1082, 2012.

\bibitem{6129407}
J.~Yang, S.~Li, and X.~Yu, ``Sliding-mode control for systems with mismatched
  uncertainties via a disturbance observer,'' \emph{IEEE Transactions on
  Industrial Electronics}, vol.~60, no.~1, pp. 160--169, 2013.

\bibitem{9339876}
Z.-H. Wu, F.~Deng, B.-Z. Guo, C.~Wu, and Q.~Xiang, ``Backstepping active
  disturbance rejection control for lower triangular nonlinear systems with
  mismatched stochastic disturbances,'' \emph{IEEE Transactions on Systems,
  Man, and Cybernetics: Systems}, vol.~52, no.~4, pp. 2688--2702, 2022.

\bibitem{yang1994disturbance}
W.~Yang and M.~Tomizuka, ``Disturbance rejection through an external model for
  nonminimum phase systems,'' \emph{Journal of Dynamic Systems, Measurement,
  and Control}, vol. 116, no.~1, pp. 39--44, 03 1994.

\bibitem{10.1115/1.2896180}
M.~Tomizuka, K.~Chew, and W.~Yang, ``Disturbance rejection through an external
  model,'' \emph{Journal of Dynamic Systems, Measurement, and Control}, vol.
  112, no.~4, pp. 559--564, 12 1990.

\bibitem{d2013control}
B.~d{'}Andr{\'e}a Novel and M.~De~Lara, \emph{Control theory for
  engineers}.\hskip 1em plus 0.5em minus 0.4em\relax Springer, 2013.

\bibitem{lathi2005linear}
B.~P. Lathi and R.~A. Green, \emph{Linear systems and signals}.\hskip 1em plus
  0.5em minus 0.4em\relax Oxford University Press New York, 2005, vol.~2.

\bibitem{2002Analysis}
P.~C. Krause, O.~Wasynczuk, S.~D. Sudhoff, and S.~D. Pekarek, \emph{Analysis of
  electric machinery and drive systems}.\hskip 1em plus 0.5em minus 0.4em\relax
  John Wiley \& Sons, 2013, vol.~75.

\bibitem{2015Effective}
M.~Tuna, C.~B. Fidan, S.~Kocabey, and S.~Gorgulu, ``Effective and reliable
  speed control of permanent magnet {DC (PMDC)} motor under variable loads,''
  \emph{Journal of Electrical Engineering and Technology}, vol.~10, no.~5, pp.
  2170--2178, 2015.

\end{thebibliography}
\end{document}